\documentclass[12pt]{article}
\usepackage{amsmath,amsthm,mathrsfs}
\usepackage{amsfonts}
\usepackage{latexsym}
\usepackage{cite}

\usepackage{t1enc}
\usepackage[latin1]{inputenc}
\usepackage[english]{babel}
\usepackage[dvips]{graphicx}
\usepackage{graphicx}
\usepackage[natural]{xcolor}
\newtheorem{theorem}{Theorem}
\newtheorem{lemma}[theorem]{Lemma}

\newtheorem{conjecture}[theorem]{Conjecture}

\newcommand{\diam}{{\rm diam}}

\begin{document}

\title{On $\{1,2\}$-distance-balancedness of generalized Petersen graphs}

\author{
{\small  Gang Ma$^{a,}$\footnote{Corresponding author.\newline{\hspace*{5mm}Email addresses: math$\_$magang@163.com (G.\ Ma), jfwang@sdut.edu.cn (J.F.\ Wang), sandi.klavzar@fmf.uni-lj.si (S.\ Klav\v{z}ar).}}\;,  \ \ Jianfeng Wang$^{a}$, \ \ Sandi Klav\v{z}ar$^{b,c,d}$}\\[2mm]
\footnotesize $^a$School of Mathematics and Statistics, Shandong University of Technology, Zibo, China\\
\footnotesize $^b$Faculty of Mathematics and Physics, University of Ljubljana, Slovenia\\
\footnotesize $^c$Faculty of Natural Sciences and Mathematics, University of Maribor, Slovenia\\
\footnotesize $^d$Institute of Mathematics, Physics and Mechanics, Ljubljana, Slovenia
}

\date{}

\maketitle

\begin{abstract}
A connected graph $G$ of diameter ${\rm diam}(G) \ge \ell$  is $\ell$-distance-balanced if $|W_{xy}|=|W_{yx}|$ for every $x,y\in V(G)$ with $d_{G}(x,y)=\ell$, where $W_{xy}$ is the set of vertices of $G$ that are closer to $x$ than to $y$. It is proved that if $k\ge 3$ and $n>k(k+2)$, then the generalized Petersen graph $GP(n,k)$ is not distance-balanced and that $GP(k(k+2),k)$ is distance-balanced. This significantly improves the main result of Yang et al.\ [Electron.\ J.\ Combin.\ 16 (2009) \#N33]. It is also proved that if $k\ge 6$, where $k$ is even, and $n>\frac{5}{4}k^2+2k$, or if $k\ge 5$, where $k$ is odd, and $n>\frac{7}{4}k^2+\frac{3}{4}k$, then $GP(n,k)$ is not $2$-distance-balanced. These results  partially resolve a conjecture of  Miklavi\v{c} and \v{S}parl [Discrete Appl.\ Math.\ 244 (2018) 143--154].
\end{abstract}

\noindent
{\bf Keywords}: Distance-balanced graph; $\ell$-distance-balanced graph; Generalized Petersen graph \\

\noindent
{\bf AMS Subj.\ Class.\ (2020)}:  05C12

\section{Introduction}
\label{S:intro}

If $G = (V(G), E(G))$ is a connected graph and $x, y\in V(G)$, then the {\it distance}, $d_{G}(x, y)$, between $x$ and $y$ is the number of edges on a shortest $x,y$-path. The {\em diameter}, $\diam(G)$, of $G$ is the maximum distance between its vertices.
The set $W_{xy}$ contains the vertices that are closer to $x$ than to $y$, that is,
$$W_{xy}=\{w\in V(G):\ d_{G}(w,x) < d_{G}(w,y)\}\,.$$
Vertices $x$ and $y$ are {\em balanced} if $|W_{xy}| = |W_{yx}|$.  For an integer $\ell \in [\diam(G)] = \{1,2,\ldots, \diam(G)\}$, the graph $G$ is $\ell$-{\em distance-balanced} if each pair $x,y$ of its vertices with $d_{G}(x,y) = \ell$ is balanced.

$1$-distance-balanced were first considered by Handa~\cite{Handa:1999} in 1999. The term ``distance-balanced'' for these graphs was proposed a decade later in~\cite{Jerebic:2008}. This has prompted a widespread research into these graphs, see~\cite{Abiad:2016, Balakrishnan:2014, Balakrishnan:2009, Cabello:2011, cavaleri-2020, fernardes-2022, ghorbani-2023, Ilic:2010, Kutnar:2006, Kutnar:2009, Kutnar:2014, Miklavic:2012, YangR:2009,ali-2021,doslic-2018, kr-2021, miklavic-2021, xu-2022}.
It was Frelih who in~\cite{Frelih:2014} extended distance-balanced graphs to $\ell$-distance balanced graphs. Also these graphs have already been investigated a lot, see~\cite{MaG:2024, MaG:arXiv2023, Frelih:2018, Jerebic:2021, Miklavic:2018}.

If $n\ge 3$ and $1\le k<n/2$, then the {\em generalized Petersen graph} $GP(n,k)$ is the graph with
\begin{align*}
V(GP(n,k)) & = \{u_i:\ i\in \mathbb{Z}_n\} \cup\{v_i:\ i\in \mathbb{Z}_n\}, \\
E(GP(n,k)) & = \{u_iu_{i+1}:\ i\in \mathbb{Z}_n\} \cup\{v_iv_{i+k}:\ i\in \mathbb{Z}_n\} \cup \{u_iv_i:\ i\in \mathbb{Z}_n\}.
\end{align*}
As it turned out, in general it is difficult to determine whether a generalized Petersen graphs is $\ell$-distance-balanced for some $\ell$. Back in the seminal paper~\cite{Jerebic:2008}, the following conjecture was proposed for the case $\ell = 1$.

\begin{conjecture}\label{C:GP} {\rm \cite{Jerebic:2008}}
For any $k\ge 2$, there exists a positive integer $n_0$ such that $GP(n,k)$ is not distance-balanced for every $n\ge n_0$.
\end{conjecture}

The conjecture has been positively resolved by Yang et al.\ as follows.

\begin{theorem}\label{T:GP-notDB} {\rm \cite{YangR:2009}}
If $k\ge 2$ and $n>6k^2$, then $GP(n,k)$ is not distance-balanced.
\end{theorem}

Miklavi\v{c} and \v{S}parl~\cite{Miklavic:2018} expanded and specified Conjecture~\ref{C:GP} to $\ell$-distance-balancedness as follows.

\begin{conjecture} {\rm \cite{Miklavic:2018}}
\label{C:GP-lDB}
Let $k\ge 2$ be an integer and let
\begin{equation*}
n_k=\left\{\begin{array}{ll}
11; & k=2,\\
(k+1)^2; & k\ odd, \\
k(k+2); & k\ge 4\ even.
\end{array}\right.
\end{equation*}
Then $GP(n,k)$ is not $\ell$-distance-balanced for any $n>n_k$ and for any $1 \le \ell < \diam(GP(n,k))$. Moreover, $n_k$ is the smallest integer with this property.
\end{conjecture}

Conjecture~\ref{C:GP-lDB} has by now been confirmed for $k=2$ in~\cite{Miklavic:2018} and for $k\in \{3,4\}$ in~\cite{MaG:arXiv2023}. These results assert that if $k=2$ and $n>11$, or $k=3$ and $n>16$, or $k=4$ and $n>24$, then $GP(n,k)$ is not distance-balanced. These are significant improvements over the bound of Theorem~\ref{T:GP-notDB} for $k\in \{2,3,4\}$. In the first main result of this paper we improve the bound of Theorem~\ref{T:GP-notDB} for an arbitrary $k$, where the case $k=2$ is included for completeness.

\begin{theorem}\label{T:oneDB}
Let $n$ and $k$ be integers, where $2\le k< n/2$.
\begin{enumerate}
\item[{\rm (i)}] 
If $k\ge 3$ and $n>k(k+2)$, then $GP(n,k)$ is not distance-balanced. In addition, $GP(k(k+2),k)$ is distance-balanced.
\item[{\rm (ii)}]
If $k=2$ and $n>10$, then $GP(n,2)$ is not distance-balanced. In addition, $GP(10,2)$ is distance-balanced.
\end{enumerate}
\end{theorem}

In our second main result we deal with $2$-distance-balancedness, where the cases $k\in \{2,3,4\}$ are included for completeness.

\begin{theorem}\label{T:twoDB}
Let $n$ and $k$ be integers, where $2\le k< n/2$.
\begin{enumerate}
\item[{\rm (i)}] 
If $k\ge 6$ and $k$ is even, then $GP(n,k)$ is not $2$-distance-balanced for any $n > \frac{5}{4}k^2+2k$.
\item[{\rm (ii)}] 
If $k\ge 5$ and $k$ is odd, then $GP(n,k)$ is not $2$-distance-balanced for any $n>\frac{7}{4}k^2+\frac{3}{4}k$.
\item[{\rm (iii)}] 
If $k=2$ and $n>10$, or $k=3$ and $n>10$, or $k=4$ and $n>21$, then $GP(n,k)$ is not $2$-distance-balanced.
In addition, $GP(10,2)$, $GP(10,3)$, and $GP(21,4)$ are $2$-distance-balanced.
\end{enumerate}
\end{theorem}

Proofs of Theorems~\ref{T:oneDB} and~\ref{T:twoDB} are respectively given in Sections~\ref{S:oneDB} and~\ref{S:twoDB}.

\section{Proof of Theorem~\ref{T:oneDB}}
\label{S:oneDB}

Let $x,y$ be vertices of a graph $G$. In addition to the already defined sets $W_{xy}$ and $W_{yx}$, let
$${}_xW_y=\{w\in V(G):\ d_{G}(w,x) = d_{G}(w,y)\}\,.$$
Clearly, $|W_{xy}|+|W_{yx}|+|{}_xW_y| = |V(G)|$, which in turn implies the following simple, but useful fact.

\begin{lemma}\label{L:non-L-DB}
Let $x,y$ be vertices of a graph $G$ with $d_{G}(x,y)=\ell$,  where $1\le\ell\le\diam(G)$.
If $2|W_{xy}|+|{}_xW_y| > |V(G)|$, then $G$ is not $\ell$-distance-balanced.
\end{lemma}

As already mentioned, Conjecture~\ref{C:GP-lDB} holds true for $k=2$. Moreover, $GP(11,2)$ is not distance-balanced, but $GP(10,2)$ is distance-balanced, see~\cite[Table~1]{Miklavic:2018}). These results cover the case $k=2$ of Theorem~\ref{T:oneDB}.

In the rest we assume that $k\ge 3$ and $n\ge k(k+2)$. We consider the vertices $u_0$ and $v_0$, and the corresponding sets $W_{u_0v_0}$, $W_{v_0u_0}$, and ${}_{u_0}W_{v_0}$.

\medskip\noindent
{\bf Case 1}: $k$ even, $k\ge 4$. In this case we have
\begin{itemize}
\item $u_i,u_{-i}\in W_{u_0v_0}$ when $0\le i\le \frac{k}{2}$; there are $2\frac{k}{2} + 1 = k + 1$ such vertices.
\item $u_i,u_{-i}\in {}_{u_0}W_{v_0}$ when $i=\frac{k+2}{2}$; there are two such vertices.
\item $u_i,u_{-i}\in W_{v_0u_0}$ when $\frac{k+2}{2}< i\le \frac{n}{2}$; there are $n - (k+3)$ such vertices.
\end{itemize}
\noindent
{\bf Subcase 1.1}: $n \bmod k =0$. In this subcase we get
\begin{itemize}
\item $v_{ik}\in W_{v_0u_0}$ when $0\le i\le\frac{n}{k}-1$; there are $\frac{n}{k}$ such vertices.
\item $\{v_i:\ 0\le i\le n-1\}\|setminus \{v_{ik}:\ 0\le i\le\frac{n}{k}-1\}\subset W_{u_0v_0}$; there are $n - \frac{n}{k}$ such vertices.
\end{itemize}
From the above we obtain
\begin{align*}
|W_{v_0u_0}|-|W_{u_0v_0}| & = \left[n-(k+3)+\frac{n}{k}\right] - \left[(k+1)+(n-\frac{n}{k})\right]\\
&=\frac{2n}{k}-2k-4.
\end{align*}

If $n>k(k+2)$, then $\frac{2n}{k}-2k-4 > 0$ and hence $|W_{v_0u_0}|>|W_{u_0v_0}|$. We can conclude
that $GP(n,k)$ is not distance-balanced if $n>k(k+2)$.

Assume now that $n=k(k+2)$. Then $\frac{2n}{k}-2k-4=0$ and hence $|W_{v_0u_0}|=|W_{u_0v_0}|$.
Since any two adjacent vertices from the set $\{u_i:\ 0\le i\le n-1\}$ as well as any two adjacent vertices from $\{v_i:\ 0\le i\le n-1\}$ are symmetrical, we can conclude that  $GP(k(k+2),k)$ is distance-balanced.

\medskip\noindent
{\bf Subcase 1.2}: $n \bmod k \neq 0$. \\
In this subcase we have $n \bmod 2k\neq 0$.  If $n>k(k+2)$, then
\begin{itemize}
\item $v_{ik},v_{-ik}\in W_{v_0u_0}$ when $0\le i\le\lfloor\frac{n}{2k}\rfloor$; there are $2\lfloor\frac{n}{2k}\rfloor + 1$ such vertices.
\end{itemize}
Hence $|W_{v_0u_0}| \ge n-(k+3)+(2\lfloor\frac{n}{2k}\rfloor+1)$ and $|{}_{u_0}W_{v_0}|\ge 2$. From this, we can estimate as follows:
\begin{align*}
2|W_{v_0u_0}|+|{}_{u_0}W_{v_0}| & \ge 2\left[n-(k+3)+(2\left\lfloor\frac{n}{2k}\right\rfloor+1)\right]+2\\
&=2n+4\left\lfloor\frac{n}{2k}\right\rfloor-2k-2\\
&\ge 2n+4\left(\frac{k+2}{2}\right)-2k-2\\
&=2n+2>2n.
\end{align*}
Applying Lemma~\ref{L:non-L-DB} we can conclude that $GP(n,k)$ is not distance-balanced.

\medskip\noindent
{\bf Case 2}: $k$ odd, $k\ge 3$. Now we obtain
\begin{itemize}
\item $u_i,u_{-i}\in W_{u_0v_0}$ when $0\le i\le \frac{k+1}{2}$; there are $2(\frac{k+1}{2}) + 1 = k+2$ such vertices.
\item $u_i,u_{-i}\in W_{v_0u_0}$ when $\frac{k+1}{2}< i\le \frac{n}{2}$; there are $n - (k+2)$ such vertices.
\end{itemize}
\noindent
{\bf Case 2.1}: $n \bmod k = 0$. In this subcase we have
\begin{itemize}
\item $v_{ik}\in W_{v_0u_0}$ when $0\le i\le\frac{n}{k}-1$; there are $\frac{n}{k}$ such vertices.
\item $\{v_i:\ 0\le i\le n-1\}\setminus \{v_{ik}:\ 0\le i\le\frac{n}{k}-1\}\subset W_{u_0v_0}$; there are $n-\frac{n}{k}$ such vertices.
\end{itemize}
By the above it follows that
\begin{align*}
|W_{v_0u_0}|-|W_{u_0v_0}| & = \left[n-(k+2)+\frac{n}{k}\right] - \left[(k+2)+(n-\frac{n}{k})\right] \\
& = \frac{2n}{k}-2k-4.
\end{align*}
If $n>k(k+2)$, then $|W_{v_0u_0}|-|W_{u_0v_0}|>0$ and $GP(n,k)$ is not distance-balanced. If $n=k(k+2)$, then $|W_{v_0u_0}|-|W_{u_0v_0}|=0$. Since any two adjacent vertices from $\{u_i:\ 0\le i\le n-1\}$ as well as any two adjacent vertices from  $\{v_i:\ 0\le i\le n-1\}$ are symmetrical, we can deduce that $GP(k(k+2),k)$ is distance-balanced.

\medskip\noindent
{\bf Case 2.2}: $n \bmod k \neq 0$. \\
Now we have $n \bmod 2k\neq 0$. Assume that $n>k(k+2)$. Then 
\begin{itemize}
\item $v_{ik},v_{-ik}\in W_{v_0u_0}$ when $0\le i\le\lfloor\frac{n}{2k}\rfloor+1$; there are $2( \lfloor\frac{n}{2k}\rfloor +1) + 1$ such vertices.
\end{itemize}
Having in mind that $k$ is odd, we have $\lfloor\frac{n}{2k}\rfloor\ge\frac{k+1}{2}$. From here we can estimate as follows:
\begin{align*}
2|W_{v_0u_0}|+|{}_{u_0}W_{v_0}| & \ge 2\left[n-(k+2)+(2\left\lfloor\frac{n}{2k}\right\rfloor+3)\right]+0\\
& = 2n+4\left\lfloor\frac{n}{2k}\right\rfloor-2k+2\\
&\ge 2n+4\left(\frac{k+1}{2}\right)-2k+2\\
&=2n+4 > 2n.
\end{align*}
Using Lemma~\ref{L:non-L-DB} once more we infer that also in this case $GP(n,k)$ is not distance-balanced. This completes the proof of Theorem~\ref{T:oneDB}.

\section{Proof of Theorem~\ref{T:twoDB}}
\label{S:twoDB}
For the case $k=2$, Theorem~\ref{T:twoDB} holds because Conjecture~\ref{C:GP-lDB} is right for $k=2$ \cite{Miklavic:2018} and
the fact that $GP(11,2)$ is not $2$-distance-balanced, but $GP(10,2)$ is $2$-distance-balanced (see Table 1 of \cite{Miklavic:2018}).
For the case $k=3$, Theorem~\ref{T:twoDB} holds because Conjecture~\ref{C:GP-lDB} is right for $k=3$ \cite{MaG:arXiv2023} and
the fact that $GP(n,3)$ is not $2$-distance-balanced when $11\le n\le 16$, but $GP(10,3)$ is $2$-distance-balanced (see Table 1 of \cite{Miklavic:2018}).
For the case $k=4$, Theorem~\ref{T:twoDB} holds because Conjecture~\ref{C:GP-lDB} is right for $k=4$ \cite{MaG:arXiv2023} and
the fact that $GP(n,4)$ is not $2$-distance-balanced when $22\le n\le 24$, but $GP(21,4)$ is $2$-distance-balanced (see Table 1 of \cite{Miklavic:2018}).


In the rest we assume that $k\ge 5$. Note that $d(u_0,v_{-k})=2$ and $v_{-k}=v_{n-k}$. We will compute $|W_{v_{-k}u_0}|$ and $|{}_{u_0}W_{v_{-k}}|$.
Two cases are discussed according to the parity of $k$.

\medskip\noindent
{\bf Case 1}: $k$ is even, $k\ge 6$, and $n > \frac{5}{4}k^2+2k$. \\
We distinguish three subcases which are separated according to which vertices are being addressed. 

\medskip\noindent
{\bf Subcase 1.1}: Vertices $u_{-i}$ and $v_{-i}$, where $1\le i\le k-1$. \\
Then $u_{-i}\in W_{u_0v_{-k}}$ and $v_{-i}\in W_{u_0v_{-k}}$ when if $1\le i\le\frac{k}{2}$, and $u_{-i}\in W_{v_{-k}u_0}$ and $v_{-i}\in {}_{u_0}W_{v_{-k}}$ when $\frac{k+2}{2}\le i\le k-1$.  So, there are $\frac{k}{2}-1$ such vertices which are in $W_{v_{-k}u_0}$ and $\frac{k}{2}-1$ such vertices which are in ${}_{u_0}W_{v_{-k}}$.

\medskip\noindent
{\bf Subcase 1.2}: Vertices $u_{i}$, where $0\le i\le n-k$. \\
For $0\le i\le k$ we have $u_{i}\in W_{u_0v_{-k}}$ when  $0\le i\le \frac{k}{2}+1$, and $u_{i}\in {}_{u_0}W_{v_{-k}}$ when $\frac{k}{2}+2\le i\le k$. Thus, there are $\frac{k}{2}-1$ such vertices which are in ${}_{u_0}W_{v_{-k}}$.

For $k+1\le i\le n-k$ we have $u_i\in {}_{u_0}W_{v_{-k}}$ or $u_i\in W_{v_{-k}u_0}$. We first consider the vertices $u_i$ such that $u_i\in W_{v_{-k}u_0}$.
Note that if $n-2k<i\le n-k$, then $u_i\in W_{v_{-k}u_0}$. Let $t$ be the largest integer such that the maximum distance of a $v_{n-k},u_i$-path is less than the minimum distance of a $u_0,u_j$-path, where $n-(t+1)k < i,j\le n-tk$. That is, $t$ is the maximal integer such that
\begin{align*}
  (t-1)+1+\frac{k}{2} < \left\lfloor\frac{n-tk}{k}\right\rfloor + 2 & \iff \\
  (t-1)+1+\frac{k}{2} < \left\lfloor\frac{n}{k}\right\rfloor-t+2 & \iff \\
  t<\frac{1}{2}\left\lfloor\frac{n}{k}\right\rfloor-\frac{k}{4}+1.
\end{align*}
Because $t$ is the largest integer satisfying the above inequality, we get 
$$t\ge\frac{1}{2}\left(\frac{n}{k}-1\right)-\frac{k}{4}+1=\frac{n}{2k}-\frac{k}{4}+\frac{1}{2}.$$
By the definition of $t$, if $1\le s\le t$, then $u_i\in W_{v_{-k}u_0}$, where $n-(s+1)k < i\le n-sk$.
That is, $u_i\in W_{v_{-k}u_0}$ for any $n-(t+1)k < i\le n-k$,  and there are $kt\ge k(\frac{n}{2k}-\frac{k}{4}+\frac{1}{2})$
such vertices which are in $W_{v_{-k}u_0}$.

Note that if $1\le j\le k$, then the difference of the distance of a $v_{n-k},u_{n-(t+1)k+j}$-path, and the distance of a $v_{n-k},u_{n-(t+2)k+j}$-path is $-1$. So, among the vertices $u_i$, where $n-(t+2)k<i\le n-(t+1)k$, 
there are at most two vertices which are not in $W_{v_{-k}u_0}$. That is, there are at least $k-2$ vertices among these which are in $W_{v_{-k}u_0}$. Using similar discussions we can get that the
number of vertices $u_i$, where $k<i\le n-(t+1)k$, which are in $W_{v_{-k}u_0}$, is at least
$$(k-2)+(k-4)+\cdots +2=\frac{k(k-2)}{4}.$$

Among the vertices $u_{i}$, where $0\le i\le n-k$, there are at least $k(\frac{n}{2k}-\frac{k}{4}+\frac{1}{2})+\frac{k(k-2)}{4}$
vertices which are in $W_{v_{-k}u_0}$,
and $n-\frac{3}{2}k-1-k(\frac{n}{2k}-\frac{k}{4}+\frac{1}{2})-\frac{k(k-2)}{4}$ vertices
which are in ${}_{u_0}W_{v_{-k}}\cup W_{v_{-k}u_0}$ and not counted in $W_{v_{-k}u_0}$.

\medskip\noindent
{\bf Subcase 1.3}:  Vertices $v_{i}$, where $0\le i\le n-k$. \\
Firstly, consider vertices $v_{sk}$ such that $v_{sk}\in {}_{u_0}W_{v_{-k}}$.
Note that $v_0\in {}_{u_0}W_{v_{-k}}$.
Let $t$ be the largest integer such that the maximum distance of a $u_0,v_{tk}$-path is less than or equal to the minimum distance of a $v_{n-k},v_{tk}$-path. That is, $t$ is the largest integer such that
$$  t+1 \le \left\lfloor\frac{n-k-tk}{k}\right\rfloor \iff
t+1\le \left\lfloor\frac{n}{k}\right\rfloor-1-t \iff 
t\le \frac{1}{2}\left\lfloor\frac{n}{k}\right\rfloor-1\,.
$$
Because $t$ is the largest integer satisfying the above inequality, we get
$$t>\frac{1}{2}\left(\frac{n}{k}-1\right)-1= \frac{n}{2k}-\frac{3}{2}.$$
By the definition of $t$ we have $v_{sk}\in {}_{u_0}W_{v_{-k}}$ if $0\le s\le t$.
That is, there are $t+1>\frac{n}{2k}-\frac{1}{2}$ such vertices which are in ${}_{u_0}W_{v_{-k}}$.

Secondly, consider vertices $v_{n-k-sk}$, such that $v_{n-k-sk}\in W_{v_{-k}u_0}$.
Note that $v_{n-k}\in W_{v_{-k}u_0}$.
Let $t$ be the largest integer such that the maximum distance of a $v_{n-k},v_{n-k-tk}$-path is less than
the minimum distance of a $u_0,v_{n-k-tk}$-path. So  $t$ is the largest integer such that
$$
t <  \left\lfloor\frac{n-k-tk}{k}\right\rfloor+1 \iff
t <  \left\lfloor\frac{n}{k}\right\rfloor-1-t+1 \iff
t < \frac{1}{2}\left\lfloor\frac{n}{k}\right\rfloor\,.
$$
Because $t$ is the largest integer satisfying the above inequality, it can be concluded that
$$t\ge \frac{1}{2}\left(\frac{n}{k}-1\right)=\frac{n}{2k}-\frac{1}{2}.$$
By the definition of $t$ we get that $v_{n-k-sk}\in W_{v_{-k}u_0}$ for $0\le s\le t$.
That is, there are $t+1\ge\frac{n}{2k}+\frac{1}{2}$ such vertices which are in $W_{v_{-k}u_0}$.

Thirdly, consider vertices $v_i$ with $0<i<n-k$, $i\neq sk$, and $i\neq n-k-sk$, such that $v_i\in {}_{u_0}W_{v_{-k}}$. Note that $v_i\in {}_{u_0}W_{v_{-k}}$ if $n-2k< i< n-k$. Let $t$ be the largest integer such that
the maximum distance of a $v_{n-k},v_i$-path is less than or equal to the minimum distance
of a $u_0,v_j$-path, where $n-(t+1)k < i,j\le n-tk$. IN other words, $t$ is the largest integer such that
\begin{align*}
(t-1)+\frac{k}{2}+2\le\left\lfloor\frac{n-tk}{k}\right\rfloor+1 & \iff \\
(t-1)+\frac{k}{2}+2\le\left\lfloor\frac{n}{k}\right\rfloor-t+1 & \iff  \\
  t\le \frac{1}{2}\left\lfloor\frac{n}{k}\right\rfloor-\frac{k}{4}.
\end{align*}
Because $t$ is the largest integer satisfying the above inequality, we can conclude that 
$$t>\frac{1}{2}\left(\frac{n}{k}-1\right)-\frac{k}{4}=\frac{n}{2k}-\frac{k}{4}-\frac{1}{2}.$$
By the definition of $t$, if $1\le s\le t$, then $v_i\in {}_{u_0}W_{v_{-k}}$, where $n-(s+1)k< i< n-sk$.
That is, there are $t(k-1)>(\frac{n}{2k}-\frac{k}{4}-\frac{1}{2})(k-1)$
such vertices which are in ${}_{u_0}W_{v_{-k}}$.

If $1\le j< k$, then the difference between the distance of a $v_{n-k},v_{n-(t+1)k+j}$-path and the distance of a $v_{n-k},v_{n-(t+2)k+j}$-path is $-1$. So among the vertices $v_i$ with $n-(t+2)k<i< n-(t+1)k$,
there are at most
two vertices which are not in ${}_{u_0}W_{v_{-k}}$. That is, there are at least $k-3$ vertices among the vertices
$v_i$, where $n-(t+2)k<i< n-(t+1)k$, which are in ${}_{u_0}W_{v_{-k}}$. Similarly we can get that the
number of vertices $v_i$ ($0<i< n-(t+1)k$, where $i\neq sk$  and $i\neq n-k-sk$, which are in ${}_{u_0}W_{v_{-k}}$, is at least
$$(k-3)+(k-5)+\cdots+1=\frac{(k-2)^2}{4}.$$

Among the vertices $v_{i}$, where $0\le i\le n-k$, there are at least $\frac{n}{2k}+\frac{1}{2}$
vertices which are in $W_{v_{-k}u_0}$ and more than 
$$\left(\frac{n}{2k}-\frac{1}{2}\right)+\left[\left(\frac{n}{2k}-\frac{k}{4}-\frac{1}{2}\right)(k-1)+\frac{(k-2)^2}{4}\right]$$ 
vertices which are in ${}_{u_0}W_{v_{-k}}$.

Combining the above three subcases, we obtain that 
\begin{align*}
|W_{v_{-k}u_0}| & \ge \left(\frac{k}{2}-1\right) + \left[k\left(\frac{n}{2k}-\frac{k}{4}+\frac{1}{2}\right) + \frac{k(k-2)}{4}\right] + \left(\frac{n}{2k}+\frac{1}{2}\right)\\
&=\frac{n}{2}+\frac{n}{2k}+\frac{k}{2}-\frac{1}{2}, 
\end{align*}
which in turn implies that the number of vertices in ${}_{u_0}W_{v_{-k}}\cup W_{v_{-k}u_0}$ which are not counted in $|W_{v_{-k}u_0}|$
is at least
\begin{align*}
\left(\frac{k}{2}-1\right) & + \left[n-\frac{3}{2}k-1 - k\left(\frac{n}{2k}-\frac{k}{4}+\frac{1}{2}\right)-\frac{k(k-2)}{4}\right] \\
&\phantom{00}+\left(\frac{n}{2k}-\frac{1}{2}\right) + \left[\left(\frac{n}{2k}-\frac{k}{4}-\frac{1}{2}\right) (k-1)+\frac{(k-2)^2}{4}\right]\\
&=n-\frac{9}{4}k-1.
\end{align*}
Therefore, 
\begin{equation*}
\begin{split}
2|W_{v_{-k}u_0}|+|{}_{u_0}W_{v_{-k}}| & \ge 2\left(\frac{n}{2}+\frac{n}{2k}+\frac{k}{2}-\frac{1}{2}\right) + \left(n-\frac{9}{4}k-1\right) \\
&=2n+\frac{n}{k}-\frac{5}{4}k-2.
\end{split}
\end{equation*}
Since $n>\frac{5}{4}k^2+2k$, we get $2|W_{v_{-k}u_0}|+|{}_{u_0}W_{v_{-k}}|>2n$. Lemma~\ref{L:non-L-DB} yields that $GP(n,k)$ is not $2$-distance-balanced.

\medskip\noindent
{\bf Case 2}: $k$ is odd, $k\ge 5$, and $n>\frac{7}{4}k^2+\frac{3}{4}k$. \\
Just as in Case 1, we are going to distinguish three subcases separated according to which vertices are being addressed. 

\medskip\noindent
{\bf Subcase 2.1}: Vertices $u_{-i}$ and $v_{-i}$, where $1\le i\le k-1$. \\
If $1\le i <\frac{k+1}{2}$, then $u_{-i}\in W_{u_0v_{-k}}$ and $v_{-i}\in W_{u_0v_{-k}}$. If $i=\frac{k+1}{2}$, then $u_{-i}\in {}_{u_0}W_{v_{-k}}$ and $v_{-i}\in {}_{u_0}W_{v_{-k}}$, and thus 
there are two such vertices  in ${}_{u_0}W_{v_{-k}}$. If $\frac{k+1}{2}< i\le k-1$, then $u_{-i}\in W_{v_{-k}u_0}$ and $v_{-i}\in {}_{u_0}W_{v_{-k}}$. So,
there are $\frac{k-3}{2}$ such vertices in $W_{v_{-k}u_0}$ and $\frac{k-3}{2}$ such vertices in ${}_{u_0}W_{v_{-k}}$.

\medskip\noindent
{\bf Subcase 2.2}: Vertices $u_i$, where $0\le i\le n-k$. \\
If $0\le i\le k$, then $u_{i}\in W_{u_0v_{-k}}$ when $0\le i\le \frac{k+1}{2}$, and
$u_{i}\in {}_{u_0}W_{v_{-k}}$ when  $\frac{k+3}{2}\le i\le k$. Thus, there are $\frac{k-1}{2}$ such vertices which are in ${}_{u_0}W_{v_{-k}}$.

If $k+1\le i\le n-k$, then $u_i\in {}_{u_0}W_{v_{-k}}$ or $u_i\in W_{v_{-k}u_0}$.
We first consider the vertices $u_i$ such that $u_i\in W_{v_{-k}u_0}$.
Note that if $n-2k<i\le n-k$, then $u_i\in W_{v_{-k}u_0}$.
Let $t$ be the largest integer such that
the maximum distance of a $v_{n-k},u_i$-path is less than the minimum distance
of a $u_0,u_i$-path, where $n-(t+1)k < i\le n-tk$. In other words, $t$ is the largest integer such that
\begin{align*}
  (t-1)+1+\frac{k+1}{2}<\left\lfloor\frac{n-tk}{k}\right\rfloor+2 & \iff \\
  (t-1)+1+\frac{k+1}{2}<\left\lfloor\frac{n}{k}\right\rfloor-t+2 & \iff \\
  t<\frac{1}{2}\left\lfloor\frac{n}{k}\right\rfloor-\frac{k}{4}+\frac{3}{4}.
\end{align*}
Because $t$ is the largest integer satisfying the above inequality, we get
$$t\ge \frac{1}{2}\left(\frac{n}{k}-1\right)-\frac{k}{4}+\frac{3}{4}=\frac{n}{2k}-\frac{k}{4}+\frac{1}{4}.$$

By the definition of $t$, if $1\le s\le t$, then $u_i\in W_{v_{-k}u_0}$, where $n-(s+1)k < i\le n-sk$.
That is, $u_i\in W_{v_{-k}u_0}$ for any $n-(t+1)k < i\le n-k$,  and there are $kt\ge k(\frac{n}{2k}-\frac{k}{4}+\frac{1}{4})$
such vertices which are in $W_{v_{-k}u_0}$.

If $1\le j\le k$, then the difference between the distance of a $v_{n-k},u_{n-(t+1)k+j}$-path and the distance of a $v_{n-k},u_{n-(t+2)k+j}$-path is $-1$. Hence, among the vertices $u_i$, where $n-(t+2)k<i\le n-(t+1)k$, 
there are at most two vertices which are not in $W_{v_{-k}u_0}$. That is, there are at least $k-2$ vertices among these vertices which are in $W_{v_{-k}u_0}$. Similarly, the
number of vertices $u_i$, where $k<i\le n-(t+1)k$, which are in $W_{v_{-k}u_0}$, is at least
$$(k-2)+(k-4)+\cdots +1=\frac{(k-1)^2}{4}.$$
Among the vertices $u_{i}$, where $0\le i\le n-k$, there are at least $k(\frac{n}{2k}-\frac{k}{4}+\frac{1}{4})+\frac{(k-1)^2}{4}$
vertices which are in $W_{v_{-k}u_0}$, and 
$$n-\frac{3}{2}k-\frac{1}{2}-k\left(\frac{n}{2k}-\frac{k}{4}+\frac{1}{4}\right)-\frac{(k-1)^2}{4}$$ vertices
which are in ${}_{u_0}W_{v_{-k}}\cup W_{v_{-k}u_0}$ and not counted in $W_{v_{-k}u_0}$.

\medskip\noindent
{\bf Subcase 2.3}: Vertices $v_i$, where $0\le i\le n-k$. \\
By a similar discussion as in Case~1.3  we obtain that $v_{sk}\in {}_{u_0}W_{v_{-k}}$ if $0\le s\le t$
($t>\frac{n}{2k}-\frac{3}{2}$),
and $v_{n-k-sk}\in W_{v_{-k}u_0}$ if $0\le s\le t$ ($t\ge \frac{n}{2k}-\frac{1}{2}$).
That is, there are $t+1>\frac{n}{2k}-\frac{1}{2}$ such vertices which are in ${}_{u_0}W_{v_{-k}}$
and $t+1\ge\frac{n}{2k}+\frac{1}{2}$ such vertices which are in $W_{v_{-k}u_0}$.

We next consider vertices $v_i$, where $0<i<n-k$, $i\neq sk$, and $i\neq n-k-sk$, such that $v_i\in {}_{u_0}W_{v_{-k}}$.
If $n-2k< i< n-k$, then $v_i\in {}_{u_0}W_{v_{-k}}$.
Let $t$ be the largest integer such that
the maximum distance of a $v_{n-k},v_i$-path is less than or equal to the minimum distance
of a $u_0,v_j$-path, where $n-(t+1)k < i,j\le n-tk$. That is, $t$ is the largest integer such that
\begin{align*}
(t-1)+\frac{k+1}{2}+2 \le \left\lfloor\frac{n-tk}{k}\right\rfloor+1 & \iff \\
(t-1)+\frac{k+1}{2}+2 \le \left\lfloor\frac{n}{k}\right\rfloor-t+1 & \iff \\
t \le \frac{1}{2}\left\lfloor\frac{n}{k}\right\rfloor-\frac{k}{4}-\frac{1}{4}.
\end{align*}
As $t$ is the largest integer satisfying the above inequality, we get
$$t > \frac{1}{2}\left(\frac{n}{k}-1\right)-\frac{k}{4}-\frac{1}{4} = \frac{n}{2k}-\frac{k}{4}-\frac{3}{4}.$$
By the definition of $t$, if $1\le s\le t$, then $v_i\in {}_{u_0}W_{v_{-k}}$ where $n-(s+1)k< i< n-sk$.
That is, there are $t(k-1)>(\frac{n}{2k}-\frac{k}{4}-\frac{3}{4})(k-1)$
such vertices which are in ${}_{u_0}W_{v_{-k}}$.

If  $1\le j< k$, then the difference between the distance of a $v_{n-k},v_{n-(t+1)k+j}$-path and the distance of a $v_{n-k},v_{n-(t+2)k}+j$-path is $-1$. So among the vertices $v_i$, where $n-(t+2)k<i< n-(t+1)k$,
there are at most two vertices which are not in ${}_{u_0}W_{v_{-k}}$. Consequently, there are at least $k-3$ vertices $v_i$, where $n-(t+2)k<i< n-(t+1)k$, which are in ${}_{u_0}W_{v_{-k}}$. Similarly, the
number of vertices $v_i$, where $0<i< n-(t+1)k$, $i\neq sk$, and $i\neq n-k-sk$, which are in ${}_{u_0}W_{v_{-k}}$, is at least
$$(k-3)+(k-5)+\cdots+2=\frac{(k-3)(k-1)}{4}.$$
Among the vertices $v_{i}$, where $0\le i\le n-k$, there are at least $\frac{n}{2k}+\frac{1}{2}$
vertices which are in $W_{v_{-k}u_0}$ and more than 
$$\left(\frac{n}{2k}-\frac{1}{2}\right) + \left[ \left(\frac{n}{2k}-\frac{k}{4}-\frac{3}{4}\right)(k-1)+\frac{(k-3)(k-1)}{4}\right]$$ vertices
which are in ${}_{u_0}W_{v_{-k}}$.

Combining the above three subcases, we obtain that 
\begin{align*}
|W_{v_{-k}u_0}| 
&\ge \frac{k-3}{2}+\left[k \left(\frac{n}{2k}-\frac{k}{4}+\frac{1}{4}\right)+\frac{(k-1)^2}{4}\right] + \left(\frac{n}{2k}+\frac{1}{2}\right)\\
&= \frac{n}{2}+\frac{n}{2k}+\frac{k}{4}-\frac{3}{4}.
\end{align*}
Consequently, the number of vertices in ${}_{u_0}W_{v_{-k}}\cup W_{v_{-k}u_0}$ which are not counted in $|W_{v_{-k}u_0}|$ is at least
\begin{align*}
& \phantom{000} \frac{k+1}{2} + \left[n-\frac{3}{2}k - \frac{1}{2} - k\left( \frac{n}{2k}-\frac{k}{4}+\frac{1}{4} \right)-\frac{(k-1)^2}{4}\right]  + \\
& \phantom{000} \left( \frac{n}{2k}-\frac{1}{2} \right) + \left[ \left( \frac{n}{2k} - \frac{k}{4} - \frac{3}{4}\right) (k-1) + \frac{(k-3)(k-1)}{4}\right] \\
& = n-\frac{9}{4}k+\frac{3}{4}.
\end{align*}
Consequently, 
\begin{align*}
2|W_{v_{-k}u_0}|+|{}_{u_0}W_{v_{-k}}| & \ge 2\left(\frac{n}{2}+\frac{n}{2k}+\frac{k}{4}-\frac{3}{4}\right) + \left(n-\frac{9}{4}k+\frac{3}{4}\right) \\
& = 2n+\frac{n}{k}-\frac{7}{4}k-\frac{3}{4}.
\end{align*}
Under the assumption $n>\frac{7}{4}k^2+\frac{3}{4}k$ we get $2|W_{v_{-k}u_0}|+|{}_{u_0}W_{v_{-k}}|>2n$, hence Lemma~\ref{L:non-L-DB} yields that $GP(n,k)$ is not $2$-distance-balanced.


\section*{Acknowledgments}

This work was supported by Shandong Provincial Natural Science Foundation of China (ZR2022MA077), the research grant NSFC (12371353) of China and IC Program of Shandong Institutions of Higher Learning For Youth Innovative Talents.
Sand Klav\v{z}ar was supported was supported by the Slovenian Research Agency ARIS (research core funding P1-0297 and projects N1-0285, N1-0218).

\end{document}